\documentclass{amsart}
\usepackage{amssymb}
\usepackage{amsbsy}
\usepackage{amscd}
\usepackage{amsmath}
\usepackage{amsthm}
\usepackage{amsxtra}
\usepackage{latexsym}

\usepackage{mathdots}

\usepackage{mathrsfs}
\usepackage{upgreek}
\usepackage{wasysym}
\usepackage[all,cmtip]{xy}
\usepackage{xcolor}
\definecolor{rouge}{rgb}{0.85,0.1,.4}
\definecolor{bleu}{rgb}{0.1,0.2,0.9}
\definecolor{violet}{rgb}{0.7,0,0.8}

\usepackage[colorlinks=true,linkcolor=bleu,urlcolor=violet,citecolor=rouge]{hyperref}

\newcommand{\mf}{\mathfrak}

\newcommand{\Z}{\mathbb{Z}}
\newcommand{\C}{\mathbb{C}}

\def\g{\mathfrak{g}}
\def\q{\mathfrak{q}}

\def\z{\mathfrak{z}}

\def\reg{{\rm reg}}

\def\sl{\mathfrak{sl}}

\def\ZS{{Z}}
\def\F{\mathcal{F}}

\def\Nc{\mathscr{N}}

\def\Vc{\mathscr{V}}
\def\leq{\leqslant}

\def\ie#1{\hskip .125em ^{e}\hskip -.125em{#1}}

\DeclareMathOperator{\Spec}{Spec}
\DeclareMathOperator{\gr}{gr}

\theoremstyle{theorem}
\newtheorem{Th}{Theorem}[section]

\newtheorem{Co}[Th]{Corollary}

\newtheorem{Conj}{Conjecture}
\theoremstyle{remark}
\newtheorem{Def}[Th]{Definition}
\newtheorem{Rem}[Th]{Remark}

\newtheorem{Quest}{Question}

\title{A remark on Mishchenko-Fomenko algebras and regular sequences}
\subjclass[2010]{17B20, 14B05}
\keywords{Mishchenko-Fomenko algebra, regular sequence, nilpotent bicone}

\author{Anne Moreau}
\address{Laboratoire Paul Painlev\'{e}, CNRS U.M.R. 8524, 
   59655 Villeneuve d'Ascq Cedex,  
   FRANCE}
\email{anne.moreau@math.univ-lille1.fr}

\begin{document}

\maketitle

\begin{abstract}
In this note, we show that the free generators of the Mishchenko-Fomenko subalgebra 
of a complex reductive 
Lie algebra, constructed by the argument shift method at a regular element, 
form a regular sequence. This result was proven by Serge Ovsienko 
in the type A at a regular and semisimple element. 
Our approach is very different, and is strongly based on geometric properties 
of the nilpotent bicone. 
\end{abstract}

\section{Introduction} 
Let $\q$ be a finite-dimensional Lie algebra over the field of complex numbers $\C$. 
The symmetric algebra $S(\q) \cong \C[\q^*]$ carries a natural Poisson structure. 
Denote by $\ZS(\q)$ the Poisson center of $S(\q)$. 
Let $\xi \in\q^*$ 
and consider the {\em Mishchenko-Fomenko subalgebra} $\mathcal{F}_{\xi}(\q)$  
of $S(\q)$ constructed by the so-called   
{\em argument shift method} \cite{MF}. 
It is generated by the $\xi$-shifts of elements  
in $\ZS(\q)$, that is, $\mathcal{F}_{\xi}(\q)$ 
is generated by all the derivatives $D_\xi^{j}(p)$ for $p \in \ZS(\q)$ 
and $j \in \{0,\ldots,\deg p-1\}$, where  
$$D_\xi^{j}(p)(x) = \frac{{\rm d}^{j}}{{\rm d}t} p(x + t \xi)|_{t=0},\qquad x \in \q^*.$$ 
It is well-known that $\mathcal{F}_\xi(\q)$ is a Poisson-commutative 
subalgebra of $S(\q)$. Furthermore, 
$${\rm trdeg}(\F_\xi(\q)) \leq \frac{\dim \q +{\rm ind}\,\q}{2} =: b(\q)$$
where ${\rm ind}\,\q$ is the {\em index} of $\q$, that is, the minimal dimension 
of the stabilizers of linear forms on $\q$ for the coadjoint representation~\cite{Di1}. 
Let $\q_\reg^*$ be the set of regular elements of $\q^*$, that is, those elements whose
stabilizer in $\q$ has the minimal dimension ${\rm ind}\,\q$, 
and $\q_{\rm sing}^*:=\q^*\setminus\q_\reg^*$.

\begin{Th}[{Panyushev-Yakimova \cite{PY}}] \label{Th:PY}
Assume that the following two conditions are satisfied:
\begin{enumerate}
\item $\ZS(\q)$ contains algebraically independent homogeneous 
elements $f_1,\ldots,f_\ell$, with $\ell = {\rm ind}\,\q$, such that 
$\sum_{i=1}^\ell \deg f_i = b(\q)$,
\item the codimension of $\q_{\rm sing}^*$ in $\q^*$ is greater than or equal to 3. 
\end{enumerate}
Then for any $\xi \in \q_\reg^*$, the Mishchenko-Fomenko algebra $\F_\xi(\q)$ is a 
polynomial algebra of Krull dimension $b(\q)$, and 
it is a maximal Poisson-commutative algebra of $S(\q)$. 
\end{Th}
Theorem \ref{Th:PY} generalizes the result of Tarasov \cite{Tar} 
for semisimple Lie algebras which are known to  
satisfy the above conditions (in \cite{Tar}, the maximality is proved 
for $\xi$ regular and semisimple). 

\begin{Quest} \label{Q:main}
In the case that $\q$ satisfies the conditions (1) and (2) of Theorem \ref{Th:PY}, 
do the free generators $D_\xi^{j}(f_i)$, 
$i=1,\ldots,\ell$, $j=0,\ldots,\deg f_i -1$ of $\F_\xi(\q)$, for $\xi \in \g_{\reg}$, form a regular sequence?
\end{Quest} 
The above question is discussed for instance in \cite[Remark 3.4]{PY}. 
The motivations come from {Gelfand-Zetlin modules}  (cf.~\cite{V91,Ov}), 
and quantizations of Mishchenko-Fomenko algebras~\cite{FFT10,Ryb06}. 

In more details, if the answer to Question \ref{Q:main} is positive and 
if $\F_\xi(\q)$ admits a quantization, that is, a (maximal) commutative subalgebra 
$\hat{\F}_\xi(\q) \subset U(\q)$ such that 
$\gr \hat{\F}_\xi(\q) \cong \F_\xi(\q)$, then $U(\q)$ is free over $\hat{\F}_\xi(\q)$ 
\cite{FutornyOvsienko05}. This implies for instance 
that any $\mu \in {\rm Specm}\, \hat{\F}_\xi(\q) $ 
lifts to a simple $U(\g)$-module, i.e., there exists 
a simple $U(\g)$-module $M_\mu$, generated by $m \in M_\mu$, 
such that for all $\gamma \in \hat{\F}_\xi(\q)$, $\gamma.m = \chi_\mu(\gamma) m$, 
where $\chi_\mu \colon \hat{\F}_\xi(\q)  \to \C$ is the character corresponding to 
$\mu$ \cite{FutornyOvsienko05}. 

From now on, let $\g$ be a reductive Lie algebra with adjoint group $G$,  
and identify $\g$ with $\g^*$ through an invariant inner product $(~|~)$.  
According to a result of Chevalley, the algebra 
$S(\g)^\g = Z(\g)$ is polynomial in $\ell$ variables, 
where $S(\g)^\g=S(\g)^G$ is the subalgebra of $S(\g)$ consisting of $G$-invariant elements. 
The nilpotent cone $\Nc$ of $\g$ is by definition the subscheme of 
$\g$ defined by the augmentation ideal of $S(\g)^\g$. 
It is well-known since Kostant \cite{Ko} that $\Nc$ is a complete intersection of codimension 
$\ell$. In other words, homogeneous generators $p_1,\ldots,p_\ell$ of 
$S(\g)^\g$ form a regular 
sequence in $\g$. 

Let us fix such generators, and order them so that 
$d_1, \ldots ,d_\ell$ is an increasing sequence with $d_i$ the degree of $p_i$. 
The Mishchenko-Fomenko algebra $\F_\xi(\g)$, for $\xi \in\g$, is then generated by 
the elements $D_\xi^{j}(p_i)$ for $i=1,\ldots,\ell$ and $j = 0,\ldots,d_i-1$. 

Let $\g_\reg$ be the set of regular elements of $\g$ and set $b:=b(\g)$.  
\begin{Th} \label{Th:main} 
Assume that $\xi \in \g_\reg$. 
Then the free generators of $\F_\xi(\g)$ 
form a regular sequence. 
Namely, for $\xi \in \g_{\reg}$, the family 
$\{D_\xi^{j}(p_i) \; ; \; i=1,\ldots,\ell, \, j = 0,\ldots,d_i-1\}$  
forms a regular sequence in $\g$. 
Equivalently, the natural morphism 
\begin{align*}\begin{array}{crcl}\sigma \colon & \g &\longrightarrow& \Spec \F_\xi(\g) \cong \C^{b}, \\
&  x  & \longmapsto& (D_\xi^{j}(p_i)(x) \; ; \; i=1,\ldots,\ell, \, j = 0,\ldots,d_i-1)\end{array}
\end{align*}
induced by the inclusion of algebras $\F_\xi(\g) \subset S(\g)$ 
is faithfully flat, that is, the extension $S(\g)$ of $\F_{\xi}(\g)$ is faithfully 
flat.
\end{Th}

As mentioned in \cite[Remark 3.4]{PY}, 
the above result was proved by Ovsienko \cite{Ov} 
for $\g=\sl_n(\C)$ and $\xi$ regular and semisimple.  

Our proof is very different. It is based on  
geometric properties of the {\em nilpotent bicone} (cf.~Definition \ref{def:bicone}) 
introduced and studied in \cite{CM1}. 
We recall in Section \ref{sec:Bicone} the main results of \cite{CM1} 
on the {nilpotent bicone}.  
As a consequence we get Theorem \ref{Th:main} for $\xi$ nilpotent and regular.  
The proof of Theorem \ref{Th:main} for an arbitrary regular $\xi$ 
is completed in Section \ref{sec:Proof}. 
In Section \ref{sec:centralizers} we discuss the case 
where $\q$  
is the centralizer a nilpotent element 
of $\g$, and formulate a conjecture. 

\subsection*{Acknowledgments}
The author is very grateful to Tomoyuki Arakawa and 
Vyacheslav Futorny for submitting this problem to her attention.
She thanks Jean-Yves Charbonnel very much for his useful remarks 
about this note. Finally, she wishes to thank the anonymous referee for his careful reading and judicious comments. 

\section{Nilpotent bicone} \label{sec:Bicone}
We assume in this section that $\g$ is simple, and 
we identify $\g$ with $\g^*$ through the Killing form $(~|~)$. 

For $p$ a homogeneous element of $S(\g)^\g$, define elements 
$p^{(j)}$ of $(S(\g) \otimes_\C S(\g))^\g \cong \C[\g\times \g]^\g$ by 
\begin{align*} 
p(s x +t y)=\sum_{j=0}^{\deg p}  p^{(j)}(x,y) s^{\deg p -j} t^j ,
\end{align*}
for all $(s,t)\in \C^2$ and $(x,y) \in\g\times\g$. Thus 
for $j = 0,\ldots,\deg p$ and $(x,\xi)\in\g\times\g$, 
$$p^{(j)} (x,\xi)= \frac{1}{j!} D_\xi^{j} (p)(x).$$

\begin{Def}[\cite{CM1}] \label{def:bicone}
The {\em nilpotent bicone} $\widetilde{\Nc}$ of $\g$  
is by definition the subscheme of $\g\times \g$ defined 
by the ideal generated  by the elements $p_i^{(j)}$ for $i =1,\ldots,\ell$ and $j = 0,\ldots,d_i$, 
$$\widetilde{\Nc} = \Spec \C[\g\times \g]/(p_i^{(j)}, \, i =1,\ldots,\ell, \, j = 0,\ldots,d_i).$$
Thus a point $(x,y) \in \g\times \g$ lies in $\widetilde{\Nc}$ 
if and only if the vector span 
generated by $x$ and $y$ is contained in nilpotent cone $\Nc$. 
\end{Def}

Set 
$$\Omega := \{(x,y)\in\g\times \g \; | \; {\rm span}_\C(x,y) \setminus (0,0) \subset \g_\reg
\text{ and } \dim {\rm span}_\C(x,y) =2\}.$$ 
Denote by $\varpi_1$ and $\varpi_2$ the first and second projections 
from $\g\times \g$ to $\g$, 
$$\begin{array}{crcl}\varpi_1 \colon & \g\times \g &\longrightarrow & \g \\
&  (x,y)& \longmapsto &x,\end{array} \qquad 
 \begin{array}{crcl}\varpi_2 \colon &\g\times \g& \longrightarrow &\g \\
 & (x,y)& \longmapsto &y \end{array}.$$ 

\begin{Th}[\cite{CM1}]  \label{Th:N}
\begin{enumerate} 
\item The nilpotent bicone is a complete intersection of dimension $3 (b-\ell)$. 
\item The images by $\varpi_1$ and $\varpi_2$ of any irreducible 
component of $\widetilde{\Nc}$ are equal to $\Nc$. 
\item  The intersection $\Omega \cap \widetilde{\Nc}$ is precisely the set 
of smooth points of $\widetilde{\Nc}$, that is, the set of 
$(x,y)$ such that the differentials of the $p_i^{(j)}$'s at $(x,y)$ are linearly independent. 
\end{enumerate}
\end{Th}
Note that the scheme $\widetilde{\Nc}$ is not reduced \cite{CM1}. 
Since the algebra $\C[\g\times\g]$ is Cohen-Macaylay,  
and since the elements 
$p_i^{(j)}$ are homogenous, part (1) of Theorem~\ref{Th:N} 
implies that any subset of the set $(p_i^{(j)}, \, i =1,\ldots,\ell, \, j = 0,\ldots,d_i)$ 
forms a regular sequence in $\g \times \g$,  \cite{Mat}. 

From Theorem \ref{Th:N}, (1) and (2), we get the following. 
\begin{Co} \label{Co:bicone_fiber}
Let $e$ be a regular nilpotent element of $\g$. 
Then the fiber $\widetilde{\Nc}_e$ of the restriction to $\widetilde{\Nc}$ of $\varpi_1$ 
(resp.~$\varpi_2$) at $e$ is a complete intersection of dimension 
$b -\ell$.  
\end{Co}

\section{Proof of Theorem \ref{Th:main}} \label{sec:Proof}
For $\xi \in \g$, denote by $\mathscr{V}_\xi(\g)$ the subscheme of $\g$ 
defined by the elements $D_\xi^{j}(p_i)$, $i=1,\ldots,\ell$, 
$j=0,\ldots,d_i-1$ of $\F_\xi(\g) \subset S(\g)$. 
Since the algebra $\C[\g]$ is Cohen-Macaylay   
and since the elements 
$D_\xi^{j}(p_i)$ are homogeneous, to prove that  for $\xi\in \g_\reg$, 
the elements $D_\xi^{j}(p_i)$, $i=1,\ldots,\ell$, 
$j=0,\ldots,d_i-1$  
form a regular sequence, 
we have to prove that  for $\xi \in \g_\reg$, the scheme $\mathscr{V}_\xi(\g)$ 
is equidimensional of dimension $b-\ell$. 
Note that each irreducible component of $\mathscr{V}_\xi(\g)$ 
has at least dimension $b-\ell$. 

Let $\g_1,\ldots,\g_s$ be the simple factors of $\g$ so that 
$\g=\z\times \g_1\times \cdots\times\g_s$ with $\z$ the center of $\g$, 
and fix $\xi \in \g_\reg$. 
From 
$$\mathscr{V}_\xi(\g) \cong \mathscr{V}_\xi(\g_1) \times \cdots \times \mathscr{V}_\xi(\g_s),$$  
we can assume that $\g$ is simple. 

Corollary \ref{Co:bicone_fiber} gives Theorem \ref{Th:main} for $\xi$ 
regular and nilpotent since for such 
$\xi$, $\mathscr{V}_\xi (\g) \cong \widetilde{\Nc}_\xi$. 
It remains to generalize the statement for an arbitrary regular $\xi \in\g$. 

Let $(e,h,f)$ be a principal $\sl_2$-triple, that is, $e$ is regular nilpotent. 
Then consider the Kostant's slice 
$$\mathscr{S}_e:=e + \g^{f},$$
where $\g^{f}$ is the centralizer of $f$ in $\g$. 
This is an affine subspace of $\g$ which consists of regular elements. 
Moreover, for any regular element $\mu \in \g$, the $G$-orbit of $\mu$ 
intersects $e+\g^{f}$ at one point \cite{Ko}. 
Thus 
$$\g_\reg= G.\mathscr{S}_e.$$ 
Since $\dim \Vc_{g.\xi}(\g) =\dim \Vc_\xi (\g)$ for any $g \in G$, 
we can assume that $\xi \in \mathscr{S}_e$. 

Let $U$ be the set of $\mu \in \g$ such that $\dim \Vc_\mu (\g) =b -\ell$. 
It is an open subset of $\g$ which contains $e$  
by Corollary~\ref{Co:bicone_fiber}. 
So $U \cap  \mathscr{S}_e$ is a nonempty subset of $\mathscr{S}_e$ which contains 
$e$. Hence for any $\mu$ in a nonempty neighborhood $W$ of $e$ in $\mathscr{S}_e$, 
$\dim \Vc_\mu (\g) = b-\ell$. 
Consider the one-parameter subgroup $\rho \colon \C^* \to G$ of $G$ defined by 
$$\forall \, t \in \C^*,\qquad \rho(t).x= t^{-2} \tilde{\rho}(t) .x$$ 
where $\tilde{\rho} \colon \C^* \to G$ is the one-parameter subgroup of 
$G$ defined by $h$. 
Then $\rho$ induces a contracting $\C^*$-action on $\mathscr{S}_e$,  
meaning that 
$$\forall \, x \in \mathscr{S}_e, \qquad  \rho(t).x \in \mathscr{S}_e, \quad \rho(t).e=e 
\quad  
\text{ and }\quad \lim_{t\to 0} \rho(t).x = e.$$ 
So for some $t\in \C^*$, $\rho(t).\xi\in W$. 
But for any $t \in \C^*$, 
$$\dim \Vc_{\xi} (\g) = \dim \Vc_{\rho(t) .\xi}(\g),$$ 
whence $\dim \Vc_\xi(\g)=b-\ell$, as desired.  
\begin{Rem} 
To generalize the statement to any arbitrary regular $\xi$, 
we have used Kostant's slice. 
This can also be deduced from the construction of Borho-Kraft \cite{Borho-Kraft79} 
about deformations of $G$-orbits. 
\end{Rem} 

It remains to prove that the morphism $\sigma$ is faithfully flat for $\xi \in \g_{\reg}$. 
As $\F_\xi(\g)$ is generated by homogeneous functions, 
the fiber at $0$ of the morphism $\sigma$ has maximal dimension. 
But by what foregoes, $\sigma^{-1}(0)\cong {\Vc}_{\xi}(\g)$ has codimension 
$b$ in $\g$. 
On the other hand, by \cite[Theorem 0.1]{PY}, $\F_\xi(\g)$ is a polynomial algebra in 
$b$ variables. So $\sigma$ is an equidimensional morphism and by 
\cite[Ch.~8, Theorem 21.3]{Mat}, $\sigma$ is a flat morphism. 
In particular by \cite[Ch.~III, Exercise 9.4]{Ha}, it is an open morphism whose 
image contains $0$. 
So $\sigma$ is surjective. Hence $\sigma$ is faithfully flat, 
according to \cite[Ch.~3, Theorem 7.2]{Mat}.

\section{Centralizers of nilpotent elements}  \label{sec:centralizers}
Other interesting examples to consider come from 
the centralizers of nilpotent elements. 

Assume that $\q$ is the centralizer $\g^{e}$ of 
a nilpotent element $e$ of $\g$. 
Then the index of $\q=\g^e$ is equal 
to $\ell$ by \cite{CM2}, 
and the algebra $S(\g^{e})^{\g^{e}}$ 
is known to be polynomial 
for a large number of element $e$ (cf.~e.g.~\cite{PPY, CM3}). 

According to the main results of \cite{CM3,CM4}, 
we have a characterization of nilpotent elements $e$ 
for which $S(\g^{e})^{\g^{e}}$ is polynomial, and   
homogeneous free generators  
form a regular sequence. They are called {\em good 
elements} in \cite{CM3}. 
In more details, 
for $p \in S(\g)$, let $\ie{p}$ 
be the initial homogeneous component of its restriction to 
$\mathscr{S}_e=e+\g^{f}$, 
with $(e,h,f)$ an $\sl_2$-triple of $\g$. 
By \cite{PPY}, if $p \in S(\g)^\g$, then 
$\ie{p} \in S(\g^{e})^{\g^{e}}$. 
Consider now the following condition: 

$(\ast)$: for some homogeneous free generators 
$q_1,\ldots,q_\ell$ 
of $S(\g)^{\g}$, we have 
$$\sum_{i=1}^{\ell} \deg \ie{q}_i=b(\g^{e}).$$  
By \cite{CM3,CM4}, the condition $(\ast)$ is satisfied 
if and only if $e$ is good. 
In addition, we have the following result: 
\begin{Th}[Arakawa-Premet \cite{AP}]
Assume that the condition $(\ast)$ and the condition 
$(2)$ of 
Theorem \ref{Th:PY} are satisfied. 
Then $\F_\xi(\g^{e})$ admits 
a quantization $\hat{\F}_\xi(\g^{e}) \subset U(\g)$.  
\end{Th}

The conditions $(\ast)$ and $(2)$ are satisfied 
for $e=0$, and (at least) in the following cases: 
$\g=\mathfrak{sl}_n(\C)$ and $e$ arbitrary (\cite{PPY,Y4}), 
$\g$ is simple not of type 
$E_8$ and $e$ is in the minimal nilpotent 
orbit of $\g$ (\cite{PPY,AP}).

The fact that homogeneous free generators of 
$S(\g^{e})^{\g^{e}}$ form a regular sequence 
when $\g=\mathfrak{sl}_n(\C)$ was known by \cite[Theorem 5.4]{PPY}. 
The fact that $\F_\xi(\g^{e})$ admits 
a quantization $\hat{\F}_\xi(\g^{e}) \subset U(\g)$ 
for $e=0$ comes from \cite{Ryb06,FFT10}. 

In view of the above remarks, we formulate a conjecture. 
\begin{Conj} \label{conjecture1}
Assume that the condition $(\ast)$ and 
the condition $(2)$ of 
Theorem~\ref{Th:PY} are satisfied. 
Then the free generators of $\F_\xi(\g^{e})$ 
form a regular sequence for any $\xi \in (\g^{e})^*_{\reg}$. 
\end{Conj}

Conjecture \ref{conjecture1} 
holds for $e=0$ (Theorem \ref{Th:main}), for regular nilpotent $e$ 
(since $\g^{e}$ is commutative in this case),  
for subregular nilpotent $e \in \mf{sl}_n(\C)$ (easy 
computations), was proved 
by Tomoyuki Arakawa and Vyacheslav Futorny for minimal nilpotent 
$e \in \mf{sl}_n(\C)$ (private communication) 
and by Wilson Fernando Mutis Cantero 
for any nilpotent $e\in \mf{gl}_n(\C)$, $n=3,4$ 
\cite{PhD:Cantero}. 

Note that $S(\g^{e})^{\g^{e}}$ is not 
always polynomial, cf.~\cite{Y3,CM3,Y4}. 
Also, even when $S(\g^{e})^{\g^{e}}$ is free, it may happen 
that the free generators 
do not form a regular sequence (cf.~\cite[Examples 7.5 and 7.6]{CM3}). 
At last, 
the codimension of $(\g^{e})_{\rm sing}^*$ in $(\g^{e})^*$ 
is not always greater than or equal to 2 (cf.~\cite{PPY}), 
even if $e$ is good \cite[Remark 7.7]{CM3}.

\end{document}